\documentclass[reqno]{article}
\usepackage{amsmath, amsfonts, amscd, eucal}
\usepackage{latexsym}
\usepackage{amssymb}
\newcommand{\be}{\begin{equation}}
      \newcommand{\ee}{\end{equation}}
      \newcommand{\ba}{\begin{eqnarray}}
       \newcommand{\ea}{\end{eqnarray}}
\newcommand{\ban}{\begin{eqnarray*}}
\newcommand{\ean}{\end{eqnarray*}}

\newcommand{\pt}{\partial}

\newcommand{\ra}{\rightarrow}

 \renewcommand{\o}[2]{\frac{#1}{#2}}
\newcommand{\hf}{\o{1}{2}}

 \newcommand{\qed}{\hspace*{\fill}\rule{3mm}{3mm}\quad \vspace{.2cm}}
 \newcommand{\Pf}{\noindent {\bf Proof:} }
 \newcommand{\Rk}{\noindent {\bf Remark} }

\newcommand{\sect}[1]{\section{#1} \setcounter{equation}{0}}

\newtheorem{theo}{Theorem}[section]

\begin{document}
\newtheorem{defn}[theo]{Definition}
\newtheorem{ques}[theo]{Question}
\newtheorem{lem}[theo]{Lemma}
\newtheorem{prop}[theo]{Proposition}
\newtheorem{coro}[theo]{Corollary}
\newtheorem{ex}[theo]{Example}
\newtheorem{note}[theo]{Note}

\title{Hitchin-Thorpe Inequality for Noncompact Einstein $4$-Manifolds}
\author{Xianzhe Dai\thanks {Math Dept, UCSB, Santa Barbara, CA 93106 \tt{Email:
dai@math.ucsb.edu}. Partially supported by NSF Grant \# DMS-0405890}  \and
Guofang Wei\thanks {Math Dept. UCSB. \tt{Email:
wei@math.ucsb.edu}. Partially supported by NSF Grant \#
DMS-0505733.}} \maketitle

\begin{abstract}  We prove a Hitchin-Thorpe inequality for noncompact
Einstein $4$-manifolds with specified asymptotic geometry at
infinity. The asymptotic geometry at
infinity is either a cusp bundle over a compact space (the fibered cusps) or a fiber bundle
over a cone with a compact fiber (the fibered boundary). Many noncompact Einstein manifolds
come with such a geometry at infinity.
\end{abstract}

\newcommand{\End}{{\rm End}}
\newcommand{\sign}{{\rm sign}}
\newcommand{\vol}{{\rm vol}}
\newcommand{\Tr}{{\rm Tr}}

\newcommand{\Ric}{\mbox{Ric}}
\newcommand{\Iso}{\mbox{Iso}}
\newcommand{\Hess}{\mbox{Hess}}
\newcommand{\divg}{\mbox{div}}
\newcommand\grd{\nabla}
\newcommand{\Dirac}{\mathcal{D}}
\newcommand{\M}{\mathcal{M}}
\newcommand{\lx}{L_X\raisebox{0.5ex}{$g$}}
\newcommand{\cuv}{\overset{\hspace{.5ex}\circ}{R}}
\newcommand{\Sr}{\mathcal{S}}
\newcommand{\ts}{\otimes}
\newcommand{\bM}{\bar{M}}
\newcommand{\gh}{\hat{g}}
\newcommand{\Pfa}{{\rm Pf}}

\def\operatorname#1{{\rm #1\,}}
\def\op{\operatorname}
\def\lam{\lambda}
\def\s{\sigma}
\def\ph{\phi}
\def\ps{\psi}
\def\e{\epsilon}
\def\d{\delta}
\def\D{\Delta}

\sect{Introduction}

Einstein manifolds are important both in mathematics and physics.
They are good candidates for canonical metrics on general
Riemannian manifolds and they are the vacuum solutions of
Einstein's field equation (with cosmological constant) in general
relativity. As a result, they are extensively studied (Cf.
\cite{besse}, \cite{lw}).

Besides space forms and irreducible symmetric spaces, a large
class of compact Einstein manifolds is given by the solution of
Calabi conjecture. Namely, a compact K\"ahler manifold with a
non-positive first Chern class admits a K\"ahler-Einstein metric
\cite{yau}, \cite{aubin}. In the case of positive first Chern
class, the work of \cite{Tian-Yau-1987, tian} says that $\mathbb{CP}^2 \# k
\overline{\mathbb{CP}}^2$ admits a K\"ahler-Einstein metric if
$3\le k\leq 8$. Other examples includes the so called
Page metric on $\mathbb{CP}^2 \# \overline{\mathbb{CP}}^2$ (only
Einstein) and certain principal torus bundles over
K\"ahler-Einstein manifolds \cite{wz}, see also the recent articles
\cite{Boyer-Galicki-preprint, Bohm-Wang-Ziller2004, Anderson2006} for Sasakian
Einstein metrics, compact homogenous Einstein manifolds, and Dehn surgery construction.

On the other hand, regarding the question of topological
obstructions, the obvious ones will be coming from that for the
Ricci curvature. Thus, if the Einstein constant is positive, the
manifold must be compact and the fundamental group is finite. If
the Einstein constant is zero, there are also obstructions coming
from Cheeger-Gromoll's splitting theorem \cite{cg}. Further, for
noncompact manifolds, the volume growth is at least linear
\cite{yau2} (see also \cite{cgt}).

In the case of compact Einstein $4$-manifolds, there are more
topological obstructions.  Berger \cite{berger} observed  that
a compact Einstein $4$-manifold must have non-negative Euler
number. Moreover, the Euler number is zero if and only if the
manifold is flat. This implies that, for example, $T^4 \# T^4$ and
$S^1 \times S^3$ are not Einstein.

Berger's  observation is considerably strengthened in the Hitchin-Thorpe
inequality \cite{h},  that for any compact oriented Einstein $4$-manifold $M^4$
\be \chi(M) \geq \o{3}{2} |\tau(M)|, \ee
where $\chi(M)$ denotes the Euler number of $M$, $\tau(M)$ the
signature. Furthermore, the equality holds if and only if either
$M$ is flat or the universal cover $\tilde{M}$ is $K3$.
The Hitchin-Thorpe inequality implies in particular that
$\mathbb{CP}^2 \# k \overline{\mathbb{CP}}^2$ cannot be Einstein
for $k \geq 9$, complementing very well the result of \cite{Tian-Yau-1987,tian}.

There are various extensions of the Hitchin-Thorpe inequality, see
\cite{Gromov1982, Kotschick1998, Sambusetti1998, Kotschick} among
others. The extensions can be summarized in the following
generalized Hitchin-Thorpe inequality due to Kotschick
\cite{Kotschick},  namely,  for any compact oriented Einstein $4$-manifold
$M^4$ \be \chi(M) \geq \o{3}{2} |\tau(M)| + \frac{1}{108 \pi^2}
\left( \lambda(M) \right)^4, \ee where $\lambda(M)$ is the volume
entropy. And equality occurs if and only if either $M$ is flat, or
the universal cover $\tilde{M}$ is $K3$ or hyperbolic.


In the case of noncompact manifolds, there are results of Tian-Yau
\cite{Tian-Yau, ty1, ty2} for the existence of K\"ahler-Einstein
metrics on the complements of a normal crossing divisor. There are
also many examples from general relativity. These are all of
finite topological type and moreover, most of them come with a special structure at infinity:
a fibration structure and an asymptotic geometry adapted to the fibration.
It should be pointed out however, that there exist
Ricci flat K\"ahler manifolds of infinite topological type
\cite{akl}.

In this note we prove a Hitchin-Thorpe inequality for noncompact
Einstein $4$-manifolds with specified asymptotic geometry at
infinity adapted to a fibration. Let $(M^n, g)$ be a noncompact complete Riemannian
manifold with finite topological type and $\bar{M}=M\cup \pt \bM$
its compactification. The metric $g$ is said to be asymptotic to a
fibered cusp if there is a defining function $x\in
C^{\infty}(\bar{M})$ of $\pt \bM$ and a fibration \be \label{bfs} F
\ra \pt \bM \stackrel{\pi}{\ra} B \ee of closed manifolds such that
\be \label{fcm} g \sim \o{dx^2}{x^2} + \pi^*g_B + x^2 g_F. \ee Here
$g_B$ is a metric on the base manifold $B$ and $g_F$ is a family of
metrics along the fibers. (The precise meaning of
asymptotic in (\ref{fcm}) and (\ref{fbm}) below will be discussed in
Section 3.) The coordinate change $x=e^{-r}$ transforms
the metric into the more standard looking
\[ g \sim dr^2 + \pi^*g_B + e^{-2r}  g_F. \]
Thus,  the geometry at infinity
is asymptotic to a fibration over the base $B$ with fibers given by
cusps over the original fiber $F$, hence the name `fibered cusps'.
Clearly the volume is  finite (assuming the dimension of the fiber is positive) in this case,
so if the metric is also Einstein,  the Einstein constant must be negative.
 Examples from \cite{Tian-Yau} have fibered cusp geometry at infinity.

The other asymptotic geometry we will consider is the so called fibered boundary metric:
 \be
\label{fbm} g \sim \o{dx^2}{x^4} +\o{ \pi^*g_B}{x^2} + g_F. \ee
Here one can use the coordinate change $x=\o{1}{r}$ in which the metric becomes
\[ g \sim  dr^2 + r^2 \pi^*g_B + g_F. \]
Hence the geometry at infinity is asymptotic to a fibration with the
original fibers $F$, but now the base is the infinite end of the
cone over the original base $B$. In this case the volume is infinite
(assuming the dimension of the base is positive) and thus the
Einstein constant could be zero or negative. The examples from
general relativity, like the Euclidean Schwarzschild metric on
$\mathbb{R}^2 \times S^2$, the Taub-NUT metric on $\mathbb{R}^4$, or
the general Gibbons-Hawking multi-center
metrics,
all have fibered boundary metric
with base $S^2$ and fiber $S^1$. The examples from \cite{ty1, ty2}
have fibration structure ($S^1$ over a smooth divisor) but the
metric is not precisely of the type we consider here.

\begin{theo} \label{mrfcm} Let $(M^4, g)$ be a noncompact complete Einstein
manifold which is asymptotic to a fibered cusp or a fibered boundary
at infinity. In the fibered boundary case, we also assume that $\dim
F >0$ (that is, we exclude the case when $F$ is a point; see below
for a separate discussion). Then
\[ \chi(M) \geq \o{3}{2}|\tau(M) + \hf\mbox{a-}\lim \eta|, \]
where $\mbox{a-}\lim \eta$ is  the adiabatic limit of eta
invariant of $\pt \bar{M}$ (for the signature operator). Moreover, the
equality holds iff $(M, g)$ is a complete Calabi-Yau manifold.
\end{theo}

\Rk One can also state an inequality with a volume entropy term. However, unlike the compact case,
it is unclear if the volume entropy here is a topological invariant.

The adiabatic limit of eta invariant of $\pt \bar{M}$ (see, e.g.  \cite{d})
encodes geometric and topological
information of the boundary fibration (at infinity). In the case when the fibration is a circle bundle over a surface,
it is given in terms of the Euler number of the circle bundle.  The case of surface bundle over a circle is more complicated.
For a torus bundle over a circle (solvmanifold) the adiabatic limit is given by certain $L$-function \cite{c}.

\begin{coro} \label{mrfbm} Let $(M^4, g)$ be a noncompact complete Einstein
manifold which is asymptotic to a fibered cusp/boundary at infinity,
with the fibration given by a circle bundle over a surface. Then
\[ \chi(M) \geq \o{3}{2}|\tau(M) - \o{1}{3} e + \sign\, e|, \]
where $e$ is the Euler number of the circle bundle. Moreover, the
equality holds iff $(M, g)$ is a complete Calabi-Yau manifold.
\end{coro}

In particular, if $M^4$ is the Taub-NUT manifold, then
$M\#(S^1\times N)$, for any closed $3$-manifold $N$, does not admit
Einstein metric with the same asymptotic geometry. Similarly, if
$M^4$ is the Taub-NUT manifold or one of the K\"ahler-Einstein
manifolds constructed in \cite{Tian-Yau}, the blowups
$M\#k\overline{\mathbb{CP}}^2$ does not admit Einstein metric with
the same asymptotic geometry for $k$ sufficiently large.

We now look at the case of fibered boundary metrics when the fiber is a single point.
In this case $B=\pt \bM$ and the geometry at infinity is asymptotically conical. That is
\[ g \sim  dr^2 + r^2 g_{\pt \bM}, \]
where $r$ can be thought as the distance from a base. Since $r$ can only change by adding a constant,
$g_{\pt \bM}$ is uniquely determined.

\begin{theo} \label{act} Let $(M^4, g)$ be a complete Einstein four manifold which is asymptotic to
a cone over $(\pt \bM, g_{\pt \bM})$. Then
\[ \chi(M) \geq  \o{1}{2\pi^2} \vol(\pt \bM) + \o{3}{2}|\tau(M) + \hf \eta(\pt \bM) | +\alpha(\pt\bM), \]
where $\eta(\pt \bM)$ is the eta invariant of $(\pt \bM, g_{\pt \bM})$ and $\alpha(\pt\bM)$ a geometric invariant
defined by
\[ \alpha(\pt\bM)= \o{1}{8\pi^2} \int_{\pt \bM} \epsilon_{abc} \omega^a\wedge [\Omega^b_c -\omega^b\wedge \omega^c] =
 \o{1}{8\pi^2} \int_{\pt \bM} \epsilon_{abc} \omega^a\wedge \Omega^b_c - \o{3}{4\pi^2} \vol(\pt \bM) \]
with $\omega^a$ denoting the dual $1$-forms of an orthonormal basis for $\pt\bM$ and $\Omega^b_c$
the $2$-form components of the curvature of  $\pt\bM$ with respect to the orthonromal basis. Moreover, the equality holds if
and only if $M$ is an asymptotically conical Calabi-Yau manifold.
\end{theo}
Note that $\alpha(S^3/\Gamma)=0$. This generalizes the previous work for ALE spaces \cite{Nakajima1990}. See also the discussion below.
\newline

One can roughly classify noncompact Einstein manifolds by their
volume growth. There are previous work concerning big volume growth.
For asymptotic locally Euclidean (hence with Euclidean volume
growth) Ricci flat $4$-manifolds with end $S^3/\Gamma$, it is proved
in \cite{Nakajima1990} that
\[ \chi(M) \geq \frac{1}{|\Gamma|} + \o{3}{2}|\tau(M) + \eta_S (S^3/\Gamma)|, \]
where $\eta_S (S^3/\Gamma)$ is the eta invariant of $S^3/\Gamma$.
Note that this class corresponds to our situation of
fibered boundary case, with the trivial fiber $F$ a single point and $B=S^3/\Gamma$.

For negative Einstein constant there are works
\cite{h2,a} on conformally compact Einstein
$4$-manifolds (hence with exponential volume growth). In this case, Anderson shows that
\[ \chi(M)-\o{3}{4\pi^2} V \geq \o{3}{2}|\tau(M)-\eta|, \]
where $V$ is the so called renormalized volume (Cf. \cite{graham})
and $\eta$ denotes the eta invariant of the conformal infinity.

Theorem~\ref{mrfcm} corresponds to the finite volume or
sub-Euclidean volume growth, while Theorem~\ref{act} corresponds to
the Euclidean volume growth.

In the process of writing this paper we learned that in the finite
volume case Yugang Zhang \cite{Zhang}  proved a similar result when
the boundary admits an injective $F$-structure and the total space
has bounded covering geometry.  While there are overlaps between
the finite volume case in Corollary~\ref{mrfbm}  and his result, as any $S^1$
bundles over a surface  have an injective $F$-structure iff  the fundamental group of the
total space is infinite, our result does cover the case of  finite fundamental group.
Our result in the infinite volume case is completely different from the corresponding case
of \cite{Zhang}.

The essential part of our proof is to extend the Gauss-Bonnet-Chern
and Hirzebruch signature formulas to complete manifolds with fibered
geometry at infinity. The index formulas we prove (Theorem
\ref{apsffcm} and Theorem \ref{apsffbm}) hold in any dimension and
should be of independent interest. Our approach is based on
application of Atiyah-Patodi-Singer index formula \cite{aps}. We use
the asymptotic structure to approximate $M$ by compact manifolds
with boundary. The boundary will in general not be totally geodesic.
Therefore, there are Chern-Simons correction terms coming from the
boundary, and analyzing these Chern-Simons correction terms consists
of the main part of the proof.

The paper is organized as follows. In Section 2 we review  the
Hitchin-Thorpe inequality for closed manifolds and the Chern-Simons
correction terms from the boundary. In Section 3, we analyze the
Chern-Simons correction term in the fibered cusp case, and fibered
boundary case and show that they limit to zero. We found out that
the language of rescaled tangent bundle introduced by Melrose
\cite{m} (see also \cite{v}) is very useful in this analysis. We
devote Section 4 to the analysis of the Chern-Simons term in the
fibered boundary case without the dimensional restriction.  Section
5 reviews the results for adiabatic limit of eta invariant.

{\it Acknowledgement}: The authors are grateful to Rafe Mazzeo, Xiaochun Rong, Gang Tian and Damin Wu
for very interesting discussions.

\sect{Chern-Simons correction term to APS}

The original Hitchin-Thorpe inequality is a beautiful application of
the Gauss-Bonnet-Chern formula and Hirzebruch's signature formula,
two special cases of the Atiyah-Singer index theorem. For a closed
oriented manifold $M$ of even dimension $n$, the Gauss-Bonnet-Chern
formula says that
\[ \chi(M) =(-1)^{n/2}\int_M \Pfa(\o{\Omega}{2\pi}), \]
where $\Omega$ is the curvature form of a Riemannian metric and
$\Pfa$ denotes the Pfaffian. For $n=4$, this gives the following
explicit formula: \ban \chi(M) & = & \o{1}{32\pi^2} \int_M
\epsilon_{abcd} \Omega_{ab}\wedge \Omega_{cd} \\
 & = & \o{1}{8\pi^2} \int_M ( |W|^2 - |Z|^2 + \o{1}{24} S^2) d\vol. \ean
 Here $W$ is the Weyl curvature, $Z$ the traceless Ricci,
 $S$ the scalar curvature, and $\epsilon_{abcd}$ denotes the totally
anti-symmetric tensor with $\epsilon_{1234}=1$ (in other words, $\epsilon_{abcd}$
is the sign of the permutation $\sigma$ where $\sigma(1)=a, \cdots, \sigma(4)=d$).

Similarly, the Hirzebruch signature formula gives
\[ \tau(M) = \int_M L(\o{\Omega}{2\pi}), \]
where $L$ denotes the $L$-polynomial. Again, in dimension $4$, the
formula simplifies to \ban  \tau(M) & = & -\o{1}{24\pi^2} \int_M
\Tr (\Omega \wedge \Omega) \\
& = & \o{1}{12\pi^2} \int_M (|W^+|^2 - |W^-|^2) d\vol .\ean

Since $|W|^2=|W^+|^2 + |W^-|^2$ and $Z=0$ for Einstein manifolds,
 the Hitchin-Thorpe inequality follows.
Furthermore, it follows that in the case of equality we must have
$S=0$, and either $W^+=0$ or $W^-=0$. That is, these must be Ricci
flat manifolds with either vanishing self dual or anti self dual
Weyl curvature. (They are shown by Hitchin \cite{h} to be either
flat or covered by $K3$.)

Assume now that $(M, g)$ is a complete noncompact manifold with
fibered geometry at infinity as defined in the previous section.
We now look at the index formula for the Euler number and
signature of such manifolds. By their topological nature, we have
\be \label{ati} \chi(M)=\chi(M_{\epsilon}), \ \ \ \tau(M)=\tau
(M_{\epsilon}), \ee for $\epsilon >0$ sufficiently small, where
$M_{\epsilon}=\{ x \geq \epsilon \}$. We are now in a position to
apply the Atiyah-Patodi-Singer index formula \cite{aps}.

If $N^n$ is an even dimensional compact oriented Riemannian
manifold with boundary $\pt N$, whose metric is the product type
near the boundary, then
\[ \chi(N) =(-1)^{n/2}\int_N \Pfa(\o{\Omega}{2\pi}), \]
and
\[ \tau(N) = \int_N L(\o{\Omega}{2\pi}) - \hf \eta(\pt N) , \]
with $\eta(\pt N)$ denoting the eta invariant of the signature
operator $A$ on the boundary with respect to the induced metric.
However, $M_{\epsilon}$ does not have product metric near its
boundary. Hence there will be Chern-Simons terms coming out as
well.

Let $P$ be an invariant polynomial of a Lie group $G$, of degree
$k$. By the Chern-Weil theory, for any $G$-connection $\omega$
with curvature $\Omega$,
\[ P(\Omega) \]
defines a characteristic form. If $\omega'$ is another
$G$-connection whose curvature form is denoted by $\Omega'$, then
their corresponding characteristic forms differ by an exact form:
\be \label{docfoc} P(\Omega') - P(\Omega) = dQ, \ee where \be
\label{csf} Q(\omega', \omega)=k \int_0^1 P(\omega'-\omega,
\Omega_t, \cdots, \Omega_t) dt . \ee Here we have denoted by
$\Omega_t$ the curvature form of the connection
$\omega_t=t\omega'+ (1-t)\omega$ interpolating between the two
connections.

Now, suppose $N$ is an compact oriented manifold with boundary whose
metric $g$ may not be product near the boundary. Then, near the
boundary $\pt N$, \[ g=dr^2 + h(r), \] where $r$ is the geodesic
distance from the boundary and $h(r)$ is the restriction of $g$ on
the constant $r$ hypersurface which is diffeomorphic to $\pt N$,
for $r$ sufficiently small. Let $g_0$ be a metric on $N$ which is
equal to $g$ except near the boundary, and is a product
sufficiently close to the boundary: \[ g_0 = dr^2 + h(0). \]
Denote by $\omega$ and $\omega_0$ the connection $1$-forms of the
Levi-Civita connections of $g$ and $g_0$, respectively. Then, by
(\ref{docfoc}),
\[ \int_N P(\Omega) - \int_N P(\Omega_0) = \int_{\pt N} Q(\omega,
\omega_0), \] where \be \label{cstfb} Q(\omega, \omega_0)=k \int_0^1
P(\theta, \Omega_t, \cdots, \Omega_t) dt, \ee and \[ \theta= \omega
- \omega_0
\] is the second fundamental form at the boundary. This is the
general form of the Chern-Simons correction to the
Atiyah-Patdi-Singer index formula for a non-product type metric. Namely,
\[ \chi(N) =(-1)^{n/2}\int_N \Pfa(\o{\Omega}{2\pi})- \int_{\pt N} Q(\omega,
\omega_0) , \]
and
\[ \tau(N) = \int_N L(\o{\Omega}{2\pi}) -  \int_{\pt N} Q(\omega,
\omega_0)- \hf \eta(\pt N). \]
Here $Q$ is associated to the Pfaffian and the L-polynomial respectively.
These formula are obtained by applying the APS index theorem to $g_0$ and
then replacing the characteristic integral of $g_0$ by that of $g$.

In dimension $4$, the Chern-Simons correction terms can be made more
explicit \cite{ch}, \cite{egh}. When $P$ is the Pfaffian, one has
\be \label{cscff} \int_{\pt N} Q(\omega, \omega_0)=
\frac{1}{32\pi^2} \int_{\pt N} \epsilon_{abcd} (2\theta^a_b\wedge
\Omega^c_d - \frac{4}{3} \theta^a_b\wedge \theta^c_e\wedge
\theta^e_d). \ee  For $P=\frac{1}{3} p_1$, it is given by \be
\label{cscfp} -\frac{1}{24\pi^2} \int_{\pt N} \Tr(\theta\wedge
\Omega) = -\frac{1}{12\pi^2} \int_{\pt N} \theta^0_i \wedge
\Omega^i_0. \ee

In the following section we study these Chern-Simons correction
terms for manifolds with fibered geometry at infinity.

\sect{Fibered geometry at infinity }

Now let $(M^n, g)$ be a noncompact complete Riemannian
manifold with finite topological type and $\bar{M}=M\cup \pt \bM$
its compactification. Moreover, there is  a fibration structure on the boundary
(at infinity)
 \[ F  \ra \pt \bM \stackrel{\pi}{\ra} B \]
with $B$, $F$ closed manifolds, as in (\ref{bfs}). Let $x$ be a boundary defining
function, i.e., $x \in  C^{\infty}(\bar{M})$, $x>0$ in $M$ and $x=0$ on $\pt \bM$; in addition
$dx$ is nowhere vanishing on $\pt \bM$.
Associated to the compactification $\bM$ of the manifold $M$ with
fibered structure at infinity (and the defining function), there is a Lie algebra of vector
fields
\[ ^{\phi}\mathcal V(\bM)= \{ \mbox{vector \ field\ } X \ \mbox{on}
\ \bM \ \mbox{tangent \ to \ the \ fibers \ at \ the \ boundary, \
and} \ X(x)=O(x^2) \}.
\]
It defines a vector bundle $^\phi T\bM$, the rescaled tangent bundle,  on $\bM$ via
\[ ^{\phi}\mathcal V(\bM)=\Gamma(^\phi T\bM). \]
If $y$, $z$ are local coordinates for the base $B$ and fiber $F$
respectively, a local frame near $\pt \bM$ for $^\phi T\bM$ is
then given by $x^2 \pt_x, \ x\pt_y, \ \pt_z$. Thus, on
$M$, where $x>0$, $^\phi T\bM$ is (non-canonically)
isomorphic to $T\bM$ (or $TM$). In turn, this induces a
non-canonical identification \[ \End(^\phi T\bM)|_{M} \cong \End (TM)
\] where different identifications differ by the adjoint action, i.e., by conjugation.
This implies that invariant polynomials are canonically
identified. For example, the trace functionals are canonically
identified: \be \label{ciot1}
\begin{array}{cccc} \Tr: & \ \End(^\phi T\bM)|_{M} & \ra &
\mathbb R \\
\parallel & \downarrow \cong & & \parallel \\
\Tr: & \ \End( TM) & \ra & \mathbb R. \end{array} \ee

A metric $g_1$ is said to be a
fibered boundary metric if there is a defining function $x\in
C^{\infty}(\bar{M})$ of $\pt \bM$
such that
\be \label{efbm} g_1= \o{dx^2}{x^4} +\o{ \pi^*g_B}{x^2} +  g_F,  \ee
where
$g_B$ is a metric on the base manifold $B$ and $g_F$ is a family of
metrics along the fibers.  Note that $g_1$ in fact defines a smooth metric on
the rescaled tangent bundle $^\phi T\bM$.

\begin{defn} \label{afbm} A metric $g$ is asymptotic to a fibered boundary metric if
\[ g=g_1 + a, \]
where $g_1$ is a fibered boundary metric defined by (\ref{efbm}) and
$a=xa_1$ where $a_1$ is a smooth section of ${\mathcal S}^2(^\phi
T\bM)$ such that $a_1(x^2 \pt_x,\  \cdot)\equiv 0$. Here ${\mathcal
S}^2$ denotes the space of symmetric two tensors.
\end{defn}

\Rk: The condition on the perturbation term in Defintion \ref{afbm} means that $a$ contains
no terms in $dx$. Thus, the normal direction to $\pt \bM$ is still given by $\pt_x$. This condition,
however, can be relaxed, see \cite{v}.

A special example of asymptotically fibered boundary metric is a metric of the form
\[ g= \o{dx^2}{x^4} +\o{ \pi^*g_B(x)}{x^2} +  g_F(x), \]
where $g_B(x)$ is a family of metrics on the base manifold $B$
depending smoothly on $x$ and $g_F(x)$ is a family of metrics along
the fibers, also depending smoothly on $x$. Many examples appear in
this form. For example, the Euclidean Schwarzschild metric on
$\mathbb{R}^2 \times S^2$, the Taub-NUT metric on $\mathbb{R}^4$, or
the general Gibbons-Hawking multi-center metrics are in this form
with fiber $S^1$.

The vector bundle $^\phi T\bM$ captures geometric information
about fibered boundary metric. The following is proved in
\cite{v}.

\begin{prop} \label{gofbm} The Levi-Civita connection for a metric
asymptotic to the fibered boundary metric is a true connection,
i.e.,
\[ \nabla^{\phi}: \ \Gamma(^\phi T\bM) \ra \Gamma(T^*\bM \otimes\,
^\phi T\bM). \] Moreover, \[ R^{\phi} \in \Gamma(\Lambda^2T^*\bM
\otimes \End(^\phi T\bM)). \]
\end{prop}

The asymptotic fibered cusp metric $g_d$ and asymptotic fibered
boundary metric $g_{\phi}$ are related by a conformal rescaling:
\[ g_d=x^2g_{\phi}, \]
and we will use this as the definition of asymptotic fibered cusp metric.
Let $^dT\bM=x^{-1}\,  ^\phi T\bM$, i.e., a local frame near the
boundary for $^dT\bM$ will be $x\pt_x, \ \pt_y, \ x^{-1}\pt_z$.
Then one also has canonical identification of the invariant
polynomials such as the trace functionals \be \label{ciot2}
\begin{array}{cccc} \Tr: & \ \End(^d T\bM) & \ra
& \mathbb R \\
\parallel & \downarrow \cong & & \parallel \\
\Tr: & \ \End( T\bM) & \ra & \mathbb R. \end{array} \ee
Furthermore, one has similarly (\cite{v})

\begin{prop} \label{gofcm} The Levi-Civita connection for a metric
asymptotic to the fibered cusp metric is a true connection, i.e.,
\[ \nabla^d: \ \Gamma(^d T\bM) \ra \Gamma(T^*\bM \otimes ^d
T\bM). \] Moreover, \[ R^d \in \Gamma(\Lambda^2T^*\bM \otimes
\End(^d T\bM)). \]
\end{prop}

Important to our consideration is the following lemma from
\cite{v} regarding the second fundamental form of asymptotic
fibered cusp metric.

\begin{lem} \label{vsffffcm} For any $T\in \Gamma( ^bT\bM)$, and $A\in
\Gamma(^dT\bM)$, \[ \nabla^d_T \o{dx}{x} (A) |_{\pt \bM} =0 .\]
\end{lem}

We are now in position to prove the following
\begin{theo} \label{apsffcm} Let $(M, g)$ be an even dimensional complete manifold
which is asymptotic to a fibered cusp metric at infinity. Then
\[ \chi(M)= (-1)^{n/2}\int_M \Pfa(\o{\Omega}{2\pi}), \]
and
\[ \tau(M) = \int_M L(\o{\Omega}{2\pi}) - \hf \mbox{a-}\lim \eta , \]
where a-$\lim \eta=\lim_{\epsilon \ra 0} \eta(\pt M_{\epsilon})$
denotes the adiabatic limit of the eta invariant.
\end{theo}

\Pf Since the proofs of both formula are similar, we do it for the
signature formula here. Applying the Atiyah-Patodi-Singer formula
with Chern-Simons correction term to $M_{\epsilon}=\{x \geq \epsilon \}$ (and $\epsilon$
sufficiently small), we have
\[ \tau(M_{\epsilon}) = \int_{M_{\epsilon}} L(\o{\Omega}{2\pi}) -\int_{\pt M_{\epsilon}}
Q - \hf \eta(\pt M_{\epsilon}) , \] where $Q$ is the Chern-Simons
terms involving the second fundamental form of $\pt M_{\epsilon}$.
By Proposition \ref{gofcm} and the discussion preceding it, we can
take $\Omega \in \Gamma(\Lambda^2T^*\bM \otimes \End(^d T\bM))$.
It follows that the first term on the right hand side of the APS
index formula has a finite limit as $\epsilon$ goes to zero. The
metric on $\pt M_{\epsilon}$ is approaching
\[ \pi^{*}g_{B} + \epsilon^{2}g_{F}= \epsilon^{2} (\epsilon^{-2}\pi^{*}g_{B} +
g_{F}). \] By the scale invariance of the eta invariance,
\[ \lim_{\epsilon \ra 0} \eta(\pt M_{\epsilon})= \mbox{a-}\lim
\eta\] is the adiabatic limit. On the other hand, the limit as
$\epsilon$ goes to zero of the Chern-Simons term is zero, since the
limit of the second fundamental form is zero as follows from Lemma
\ref{vsffffcm}. Our result follows. \qed

For fibered boundary geometry at infinity, the analysis of the
Chern-Simons term is more complicated. We will restrict ourself to
dimension $4$ in this section and leave the general discussion to
the next section. As we see from (\ref{cscff}) and (\ref{cscfp}),
this involves computing the second fundamental form $\theta$ and the
curvature form $\Omega$. Taking a cue from our treatment in the
fibered cusp case, we express both $\theta$, $\Omega$ as matrices
with respect to an orthonormal basis, but with entries differential
forms that are smooth up to the boundary at infinity $x=0$. First,
assume that $g=g_1$ is a fibered boundary metric as defined by
(\ref{efbm}). Fix a local orthonormal frame $e_a, \ e_i$ of $\pt
\bar{M}$ compactible with the submersion metric $\pi^*(g_B) + g_F$
and let $\theta^a, \theta^i$ be the dual $1$-forms, where $a$ ranges
over the coordinates of $B$ and $i$ that of $F$. Then near infinity,
\[ x^2 \pt_x, \ x\,e_a, \ e_i \]
form an orthonormal basis for the metric $g$. Computing with respect
to this basis, we find at $x=0$
\[ \theta^0_a=\theta^a, \ \ \ \theta^0_i=0. \]
Similarly, we find
\[ \Omega^a_0= f^a_{bc}(x)\theta^b \wedge \theta^c + O(x^2), \]
where $f^a_{bc}(x)=O(1)$ as $x \ra 0$. This shows that, for fibered
boundary geometry at infinity where the fibration has positive
dimensional fiber, the Chern-Simons term (Cf. (\ref{cscfp})) for the
signature vanishes in dimension $4$:
\[ \theta^0_a \wedge \Omega^a_0 = O(x^2). \]

For the Chern-Simons terms for the Euler number (\ref{cscff}), one term
involves only the second fundamental form:
\[ \epsilon_{abcd} \theta^a_b\wedge \theta^c_e \wedge \theta^e_d. \]
Since $\theta^a_b$ is zero unless one of the indices is $0$, this
term reduces to a multiple of
\[ \theta^0_1\wedge \theta^0_2 \wedge \theta^0_3 \]
which vanishes by the explicit form of $\theta^a_b$ computed above
(i.e. $\theta^0_3=0$ as $3$ is the index for the fiber coordinate
here). The other term involved is (up to a constant multiple)
\[ \epsilon_{abcd} \theta^a_b\wedge \Omega^c_d \]
which reduces to a multiple of
\[ \theta^0_1 \wedge \Omega^2_3 - \theta^0_2 \wedge \Omega^1_3. \]
Again, explicit computation gives
\[ \Omega^a_i =f^a_b(x) \theta^b \wedge \theta^3 + f^a_{bc}(x)
\theta^b \wedge \theta^c + g^a_b(x) \theta^b \wedge dx  + g(x)
\theta^3 \wedge dx,
\] where $f^a_b(x)=O(x)$. It follows then that
\[ \lim_{\epsilon \ra 0} \int_{\pt M_{\epsilon}} \theta^0_1 \wedge \Omega^2_3 - \theta^0_2 \wedge
\Omega^1_3 =0. \]

Thus, we have proved the following theorem in the special case when
$g=g_1$ is a fibered boundary metric.

\begin{theo} \label{apsffbm4} Let $(M, g)$ be a complete manifold of
dimension $4$ which is asymptotic to a fibered boundary metric at
infinity and the fiber has positive dimension. Then
\[ \chi(M)= (-1)^{n/2}\int_M \Pfa(\o{\Omega}{2\pi}), \]
and
\[ \tau(M) = \int_M L(\o{\Omega}{2\pi}) - \hf \mbox{a-}\lim \eta . \]
\end{theo}

In order to prove the theorem in general, we now consider the effect
of the perturbation term. This part of discussion is not restricted
to dimension four. Thus let $g=g_1 + a$ and denote by $\nabla$,
$\nabla^1$, the Levi-Civita connection of $g$, $g_1$ respectively.

\begin{lem} \label{pert}
Let $Q$, $Q^1$ denote the Chern-Simons correction terms with respect
to the metrics $g$, $g_1$, respectively. Then, for perturbation $a$
satisfying the condition in Definition \ref{afbm}, we have
\[ \lim_{\epsilon \rightarrow 0} \int_{\pt M_{\epsilon}} Q =
\lim_{\epsilon \rightarrow 0} \int_{\pt M_{\epsilon}} Q^1. \]
\end{lem}

\Pf Let $S=\nabla - \nabla^1$ be the difference tensor. An easy
calculation using Koszul's formula yields \ba \label{difft} g(S(X)Y,
Z) + a(\nabla^1_X Y, Z) & = & \o{1}{2} \left[ X(a(Y, Z)) + Y(a(X,
Z)) - Z(a(X, Y)) \right. \nonumber \\ & & \left. - a(X, [Y, Z]) -
a(Y, [X, Z]) + a(Z, [X, Y]) \right] , \ea for vector fields $X, Y,
Z$.

By Proposition \ref{gofbm}, $S$ is a (regular) $1$-form valued
endomorphism of $^\phi T\bM$.  Effectively, this means that in (\ref{difft})
we let $X$ be a usual vector field while letting $Y, Z$ be smooth sections
of $^\phi T\bM$, i.e., rescaled vector fields. It follows from the assumption on
the perturbation $a$ that
\[ S= xS_1 + dx\otimes S', \]
where $S_1$  is a $1$-form valued
endomorphism of $^\phi T\bM$, and $S'$ an endomorphism of $^\phi T\bM$.
The crucial point here is that the precise form of $S'$ is not important when we restrict
to $\pt M_{\epsilon}=\{x=\epsilon\}$.

Now the curvature of $g$ is related to that of $g_1$ via
\[ \Omega=\Omega_1 + [\nabla^1, S] + S^2. \]
Hence,
\[ \Omega=\Omega_1 + x \Omega' + dx \wedge \Omega'', \]
where $\Omega'$ is a (regular) $2$-form valued
endomorphism of $^\phi T\bM$ and $\Omega''$ a (regular) $1$-form valued
endomorphism of $^\phi T\bM$.

Similarly we find that the second fundamental forms of $\pt M_{\epsilon}$ with respect to the metrics $g$ and $g_1$,
respectively, differ by a term vanishing to first order of $\epsilon$:
\[ \theta=\theta_1 + \epsilon \theta', \]
where $\theta'$ is a (regular) $1$-form valued
endomorphism of $^\phi T\bM$. Our lemma follows. \qed

\sect{Chern-Simons term for fibered boundary geometry}

It turns out that Theorem \ref{apsffbm4} holds in any dimension. In
order to see this, we now discuss briefly some elementary geometry
of a fibration following \cite{bc}. Thus let $F\longrightarrow N
\stackrel{\pi}{\longrightarrow} B$ be a fibration of smooth
manifolds. It gives rise to a subbundle of $TN$, the vertical bundle
$T^V\!N$, whose section consists of vector fields of $N$ tangent to
the fibers. This leads to the exact sequence of vector bundles
\[ 0\ra T^V\!N \ra TN \ra \pi^*TB\ra 0. \]
A connection for the fibration is a splitting \be \label{cs} TN =
T^H\!N \oplus T^V\!N, \ee where $T^H\!N\cong \pi^*TB$ is the
horizontal bundle. For example, a Riemannian metric $g$ on $N$
determines such a splitting, where $T^H\!N$ is the orthogonal
complement of $T^V\!N$.

If $\nabla^F$ is a family of connections on $F$ parametrized by $B$,
then it defines a connection (still denoted by the same notation) on
$T^V\!N$ by adding
\[ \nabla^F_{X^H}Y=[X^H, Y], \]
where $X^H$ is a horizontal vector field and $Y$ vertical vector
field (a section $T^V\!N$). In particular, if $g^F$ is a family of
Riemannian metrics on $F$ (parametrized by $B$), the corresponding
Levi-Civita connections define such a connection on the vertical
bundle.

Together with a connection $\nabla^B$ (determined by a metric $g^B$
for example) on $B$, one can define a connection $\nabla$ on $TN$
which is diagonal with respect to the splitting (\ref{cs}): \be
\label{dc} \nabla=\pi^*\nabla^B \oplus \nabla^F. \ee

Now let $g^N=\pi^*g^B+g^F$ be a submersion metric on $N$. Then the
above discussion gives a diagonal connection $\nabla$ on $TN$
determined by the Levi-Civita connections of $g^B$ and $g^F$. Let
$\nabla^L$ be the Levi-Civita connection of $g^N$ and
$S=\nabla^L-\nabla$ the difference tensor. Since Levi-Civita
connections are scale invariant, the diagonal connection $\nabla$
stays the same under the adiabatic limit
$g^N_{\epsilon}=\epsilon^{-2} \pi^*g^B + g^F$.

Let $\nabla^{L, \epsilon}$ be the Levi-Civita connection of
$g^N_{\epsilon}$ and $S^{\epsilon}=\nabla^{L,\epsilon}-\nabla$ the
corresponding difference tensor. Denote by $P^H$, $P^V$ the
projections associated with the splitting (\ref{cs}). The following
observation is from \cite{bc}.

\begin{lem} For any vector field $X$ on $N$, $S(X)$ defines an odd endomorphism
of $TN$ with respect to the splitting (\ref{cs}). That is,
\[ S(X): \  T^H\!N \longrightarrow T^V\!N,  \ \ \ \  S(X): \  T^V\!N \longrightarrow T^H\!N. \]
Moreover,
\[ P^H S^{\epsilon} =\epsilon^2 P^HS, \ \ \ \  P^V S^{\epsilon} =
P^VS.\]
\end{lem}

For later purpose and also for symmetry, we paraphrase it in terms of the rescaled splitting
\be \label{rs}  ^{\epsilon}\!TN =\epsilon T^H\!N \oplus T^V\!N, \ee
and think of the connection $\nabla^{L, \epsilon}$ as a connection $^{\epsilon}\nabla^{L}$ on
 $^{\epsilon}\!TN $. (Effectively this is computing the connection with respect to an orthonormal basis
 of the adiabatic metric $g^N_{\epsilon}$ but with the crucial difference that the directional vector field is the usual
 vector field). Note that $\nabla$ stays unchanged.  Then
 \be \label{rcib}  ^{\epsilon}\nabla^{L}=\nabla + O(\epsilon). \ee

We now consider a Riemannian manifold $(M, g)$ which is asymptotic
to a fibered boundary metric at infinity. By Lemma \ref{pert}, we can actually assume that
$g$ is a fibered boundary metric. That is, we have a
fibration $F\ra \pt\bM \stackrel{\pi}{\ra} B$ and
\[ g = \frac{dx^2}{x^4} + \frac{\pi^*g^B}{x^2} + g^F. \]
Thus, near $\pt\bM$, we have a direct sum decomposition \be
\label{dffbm} ^{\phi}TM=\langle x^2 \pt_x \rangle \oplus \, x\,
T^H\!(\pt\bM) \oplus T^V\!(\pt\bM). \ee

Let $\nabla^M$ be the Levi-Civita connection of $g$. For each
hypersurface $x=\epsilon$, the metric
\[ g_0=\frac{dx^2}{\epsilon^4} + \frac{\pi^*g^B}{\epsilon^2} + g^F \]
is a product metric near $x=\epsilon$ and restricts to $x=\epsilon$
to the same metric as $g$. Let $\nabla^0$ be its Levi-Civita
connection. The difference
\[ \theta = \nabla^M - \nabla^0 \in \Omega^1(M, \End(TM)) \]
is a matrix with $1$-form entries and, when restricted to
$x=\epsilon$, has only normal components (i.e. off-diagonal with
respect to the decomposition into tangential and normal part)
determined by the second fundamental form of $x=\epsilon$. As
before, we reinterpret $\theta$ as a $1$-form taking values in
$\End(^{\phi}TM)$ and thus, $\theta$ is off-diagonal with respect to
the decomposition (\ref{dffbm}).  In fact, if we take the
orthonormal basis $x^2\pt_x, xe_i, f_{\alpha}$ where $e_i$ is (the
lift of) an orthonormal basis of $(B, g^B)$  and $ f_{\alpha}$ that
of $g^F$, then
\[  \theta_i^0= -\theta_{0}^i=\omega^i, \]
and all other components of $\theta$ vanish. Here $\omega^i$ are the
(pullback of) dual $1$-forms of $e_i$.  The crucial observation is
that $\theta$ is in block form with respect to the splitting
(\ref{dffbm}) with nontrivial entries only in the block from $
\langle x^2 \pt_x \rangle \oplus \, x\, T^H\!(\pt\bM)$ to itself.
Moreover, the nontrivial entries are (pullbacks) of forms on $B$.

\begin{lem} For a complete manifold $(M, g)$ which is asymptotic to
a fibered boundary metric at infinity and the fiber has positive dimension, and for any invariant
polynomial $P$, the Chern-Simons term vanishes at infinty:
\[ \lim_{\epsilon \ra 0} \int_{\pt M_{\epsilon}} Q =0. \]
Here $Q$ is defined in (\ref{cstfb}) and $\pt M_{\epsilon}$ is the
hypersurface $x=\epsilon$.
\end{lem}

\Pf By the discussion above, we have a similar block structure for
$\Omega_t$ with diagonal blocks from $ \langle x^2 \pt_x \rangle
\oplus \, x\, T^H\!(\pt\bM)$ to itself and from $T^V\!(\pt\bM)$ to
itself plus an error term of $O(\epsilon)$. Moreover, the diagonal
block from $ \langle x^2 \pt_x \rangle \oplus \, x\, T^H\!(\pt\bM)$
to itself involves only pullbacks of forms on $B$. It follows from
the explicit block structure of $\theta$ that
\[ P(\theta, \Omega_t, \cdots, \Omega_t)=\pi^*(\alpha) + O(\epsilon), \]
where $\alpha$ is a differential form on $B$. Hence
\[ \lim_{\epsilon \ra 0} \int_{\pt M_{\epsilon}} Q =k \int_0^1 \int_{\pt \bM}
\pi^*(\alpha)=0. \]
\qed

 Thus, we have

\begin{theo} \label{apsffbm} Let $(M, g)$ be a complete manifold
which is asymptotic to a fibered boundary metric at infinity and the fiber has positive dimension.
Then
\[ \chi(M)= (-1)^{n/2}\int_M \Pfa(\o{\Omega}{2\pi}), \]
and
\[ \tau(M) = \int_M L(\o{\Omega}{2\pi}) - \hf \mbox{a-}\lim \eta . \]
\end{theo}

\Rk When the fiber reduces to a point, one does have boundary contributions in these index formulas, which are
related to the Weyl volume of tubes invariants,  see \cite{albin1}.

 \sect{The adiabatic limit of eta invariant}

There is extensive work on the adiabatic limit of eta invariant (and
other geometric invariants) (Cf.
\cite{bc}, \cite{d} among others). In general if $M$ is a closed
oriented manifold that has a fibration structure
\begin{equation}
Y \rightarrow M\stackrel{\pi}{\rightarrow}B \label{f}
\end{equation}
and $g_{M}$ a submersion metric,
\[g_{M} = \pi^{*}g_{B} + g_{Y}, \]
then blowing up the metric in the horizontal direction by a factor
$\epsilon^{-2}$ gives us a family of metrics $g_{x}$,
\[g_{\epsilon} = \epsilon^{-2}\pi^{*}g_{B} + g_{Y}. \]
Let $A_{\epsilon}$ be the signature operator on $M$ with respect
to the adiabatic metric $g_{\epsilon}$. A general formula for
$\lim_{\epsilon\rightarrow 0}\eta(A_{\epsilon})$ is given in
\cite{d}, which, in fact, comes from a more general formula for
Dirac operators (Cf. \cite{d}). Namely, \be
\lim_{\epsilon\rightarrow 0}\eta(A_{\epsilon})=2 \int_{B} {\cal
L}(\frac{R^{B}}{2\pi})\wedge \tilde{\eta} + \eta(A_{B}\otimes \ker
A_{Y}) + 2\tau, \label{adlimg} \ee where $\tilde{\eta}$ is the the
$\tilde{\eta}$-form of Bismut-Cheeger \cite{bc}, $R^{B}$ is the
curvature tensor of $g_{B}$ and $A_{B}$ denotes the signature
operator on B and $A_{Y}$ the family of signature operators along
Y. The integer $\tau$ is a topological invariant computable from
the Leray spectral sequence.

In the case of circle bundles, i.e., $Y=S^1$, the terms on the
right hand side of (\ref{adlimg}) can be explicitly computed. For
example
\[  \tilde{\eta}= 2 ( \frac{1}{2\tanh \frac{e}{2}}-\frac{1}{e}), \]
and
\[ \tau=\sign (B_e),  \]
where $B_e$ is the quadratic form \ban
B_e: \  H^{2k-2}(B) &\otimes & H^{2k-2}(B) \ra \mathbb R \\
B_e (x \otimes y) & = & \langle xye,[B] \rangle. \ean
 Here $e$ is the Euler class of
the circle bundle. This gives us the following result of
\cite{dz}.

\begin{theo} \label{adlim}
We have \be \hf \lim_{\epsilon \ra 0} \eta (A_{\epsilon}) =
\langle L(TB) ( \frac{1}{\tanh e}-\frac{1}{e}), [B] \rangle  -
\sign (B_e). \label{adlims} \ee
\end{theo}

When $\dim B=2$, i.e., we have a circle bundle over a surface, the
formula (\ref{adlims}) gives \be \label{adlid4} \hf \lim_{\epsilon
\ra 0} \eta (A_{\epsilon})= \o{1}{3} e - \sign\, e. \ee

\sect{Proof of the theorems}

We now proceed to prove Theorem \ref{mrfcm}. By Theorems
\ref{apsffcm}, \ref{apsffbm4}, formula (\ref{adlid4}), and the decomposition of
curvature in dimension $4$, we have
\[ \chi(M)=\o{1}{8\pi^2} \int_M ( |W|^2 - |Z|^2 + \o{1}{24} S^2)
d\vol, \] and \[ \tau(M) + \o{1}{3} e - \sign\, e =\o{1}{12\pi^2}
\int_M (|W^+|^2 - |W^-|^2) d\vol  . \] The rest of the proof is
the same as in the closed case.

Note that the equality holds exactly as in the closed case, namely,
for Ricci flat manifolds with either vanishing self dual or anti
self dual Weyl curvature. Thus $M$ must be K\"ahler as follows from the same argument
of \cite{h}, and hence
Calabi-Yau.


For Theorem \ref{act}, we can no longer apply Theorem   \ref{apsffbm4}.
However, using Lemma \ref{pert},  the conformal invariant of the Pontryagin forms and the scale invariance
of the eta invariant, one still
has
\[ \tau(M) + \o{1}{2}\eta(\pt\bM) =\o{1}{12\pi^2}
\int_M (|W^+|^2 - |W^-|^2) d\vol. \]
On the other hand,
\[ \chi(M) + \lim_{\epsilon \rightarrow 0}  \int_{\pt M_{\epsilon}} Q =\o{1}{8\pi^2} \int_M ( |W|^2 - |Z|^2 + \o{1}{24} S^2)
d\vol, \]
where $Q$ is given by (\ref{cscff}):
\[  Q= \epsilon_{abcd} (2\theta^a_b\wedge
\Omega^c_d - \frac{4}{3} \theta^a_b\wedge \theta^c_e\wedge
\theta^e_d). \]
Here we emphasize that $\Omega^c_d$ denotes the two form components of the curvature of $M$.
In this case, using Lemma \ref{pert}, an explicit computation shows that
\[  \lim_{\epsilon \rightarrow 0}  \int_{\pt M_{\epsilon}} Q = \o{1}{2\pi^2} \vol(\pt\bM) +\alpha(\pt\bM). \]

\end{document}